%% file: 0404.tex
\documentclass[10pt,twoside,draft]{article}
\usepackage{epsf}
\usepackage{epic}
\usepackage{eepic}
\usepackage{amscd}
\usepackage{amssymb}
\usepackage{Latex-document}

\newtheorem{theorem}{Theorem}
\newtheorem{lemma}{Lemma}

\newtheorem{definition}{Definition}
\newcommand{\real}{{\bf R}}

\title{\bf Black Holes and the Penrose Inequality \vskip -2mm in General Relativity\vskip 6mm}
\author{Hubert L. Bray\vspace*{-0.5cm} \thanks{Mathematics Department, 2-179,
Massachusetts Institute of Technology, 77 Massachusetts Avenue,
Cambridge, MA  02139, USA. E-mail: bray@math.mit.edu}}
\date{\vspace{-8mm}}

\begin{document}

\maketitle

\thispagestyle{first} \setcounter{page}{257}

\begin{abstract}

\vskip 3mm

In a paper \cite{P} in 1973, R. Penrose made a physical argument
that the total mass of a spacetime which contains black holes with
event horizons of total area $A$ should be at least
$\sqrt{A/16\pi}$.  An important special case of this physical
statement translates into a very beautiful mathematical inequality
in Riemannian geometry known as the Riemannian Penrose inequality.
One particularly geometric aspect of this problem is the fact that
apparent horizons of black holes in this setting correspond to
minimal surfaces in Riemannian 3-manifolds. The Riemannian Penrose
inequality was first proved by G. Huisken and T. Ilmanen in 1997
for a single black hole \cite{HI} and then by the author in 1999
for any number of black holes \cite{Bray}. The two approaches use
two different geometric flow techniques.  The most general version
of the Penrose inequality is still open.

In this talk we will sketch the
author's proof by flowing Riemannian manifolds inside the class of
asymptotically flat 3-manifolds (asymptotic to $\real^3$
at infinity) which have nonnegative scalar curvature and contain minimal spheres.
This new flow of metrics has very special properties and simulates an initial
physical situation in which all of the matter falls into the black holes which
merge into a single, spherically symmetric black hole given by the Schwarzschild
metric.  Since the Schwarzschild metric gives equality in the Penrose inequality and
the flow decreases the total mass while preserving the area of the horizons of the
black holes, the Penrose inequality follows.  We will also discuss how these
techniques can be generalized in higher dimensions.

\vspace{4.5mm}

\noindent {\bf 2000 Mathematics Subject Classification:} 53, 83.

\noindent {\bf Keywords and Phrases:} Black holes, Penrose
inequality, Positive mass theorem, Quasi-local mass, General
relativity.
\end{abstract}

\section{Introduction}

\vskip -5mm \hspace{5mm}

A natural interpretation of the Penrose inequality is that the
mass contributed by a collection of black holes is (at least)
$\sqrt{A/16\pi}$, where $A$ is the total area of the event
horizons of the black holes. More generally, the question ``How
much matter is in a given region of a spacetime?'' is still very
much an open problem \cite{CY}. In this paper, we will discuss
some of the qualitative aspects of mass in general relativity,
look at examples which are informative, and sketch a proof of the
Riemannian Penrose inequality.

\subsection{Total mass in general relativity}\vskip -5mm \hspace{5mm}

Two notions of mass which are well understood in general relativity are
local energy density at a point and the total mass of an asymptotically
flat spacetime.
However, defining the mass of a region larger than a point but smaller than
the entire universe is not very well understood at all.

Suppose $(M^3,g)$ is a Riemannian 3-manifold isometrically embedded in a
(3+1) dimensional Lorentzian spacetime.  Suppose that $M^3$ has zero
second fundamental form in the spacetime.  This is a simplifying assumption
which allows us to think of $(M^3,g)$ as a ``$t=0$'' slice of the spacetime.
The Penrose inequality (which allows for $M^3$ to have general second
fundamental form) is known as the Riemannian
Penrose inequality when the second fundamental form is set to zero.

We also want to only consider $(M^3,g)$ that are asymptotically flat at
infinity, which means that for some compact set
$K$, the ``end'' $M^3 \backslash K$ is diffeomorphic to $\real^3 \backslash B_1(0)$, where
the metric $g$ is asymptotically approaching (with certain decay
conditions) the standard flat metric $\delta_{ij}$
on $\real^3$ at infinity.  The simplest example of an asymptotically
flat manifold
is $(\real^3,\delta_{ij})$ itself.  Other good examples are the conformal
metrics $(\real^3, u(x)^4 \delta_{ij})$, where $u(x)$ approaches a constant
sufficiently rapidly at infinity.  (Also, sometimes it is convenient to allow
$(M^3,g)$ to have multiple asymptotically flat ends, in which case each
connected component of $M^3 \backslash K$ must have the property described above.)

The purpose of these assumptions
on the asymptotic behavior of $(M^3,g)$ at infinity is that they
imply the existence of the limit
\begin{equation}\label{eqn:ADM_mass}
m = \frac{1}{16\pi} \lim_{\sigma\to\infty}
\int_{S_\sigma}\sum_{i,j}(g_{ij,i}\nu_j-g_{ii,j}\nu_j)\,d\mu,
\end{equation}
where $S_\sigma$ is the coordinate sphere
of radius $\sigma$, $\nu$ is the unit normal to $S_\sigma$, and $d\mu$
is the area element of $S_\sigma$ in the coordinate chart.
The quantity $m$ is called the {\bf total mass} (or ADM mass) of
$(M^3,g)$ (see \cite{ADM}, \cite{Ba2}, \cite{S}, and \cite{SY5}).

Instead of thinking of total mass as given by equation \ref{eqn:ADM_mass}, it
is better to consider the following example.
Going back to the
example $(\real^3, u(x)^4 \delta_{ij})$, if we suppose that $u(x) > 0$ has
the asymptotics at infinity
\begin {equation}\label{eqn:exp}
   u(x) = a + b/|x| + {\cal O}(1/|x|^2)
\end{equation}
(and derivatives of the ${\cal O}(1/|x|^2)$ term are ${\cal O}(1/|x|^3)$),
then the total mass of $(M^3,g)$ is
\begin{equation}\label{eqn:tmass}
   m = 2ab.
\end{equation}
Furthermore, suppose $(M^3,g)$ is any metric whose ``end'' is
isometric to $(\real^3 \backslash K,$\\$ u(x)^4 \delta_{ij})$,
where $u(x)$ is harmonic in the coordinate chart of the end
$(\real^3 \backslash K, \delta_{ij})$ and goes to a constant at
infinity. Then expanding $u(x)$ in terms of spherical harmonics
demonstrates that $u(x)$ satisfies condition \ref{eqn:exp}. We
will call these Riemannian manifolds $(M^3,g)$ {\bf harmonically
flat at infinity}, and we note that the total mass of these
manifolds is also given by equation \ref{eqn:tmass}.

A very nice lemma by Schoen and Yau is that, given any $\epsilon >0$,
it is always possible to
perturb an asymptotically flat manifold to become harmonically flat
at infinity such that the total mass changes less than $\epsilon$ and the
metric changes less than $\epsilon$ pointwise, all while maintaining
nonnegative scalar curvature (discussed in a moment).
Hence, it happens that to prove the theorems
in this paper, we only need to consider harmonically flat manifolds!
Thus, we can use equation \ref{eqn:tmass} as our definition of total mass.
As an example,
note that $(\real^3,\delta_{ij})$ has zero total mass.  Also, note
that, qualitatively, the total mass of an asymptotically flat
or harmonically flat manifold is
the $1/r$ rate at which the metric becomes flat at infinity.

\subsection{Local energy density}\vskip -5mm \hspace{5mm}

Another quantification of mass which is well understood is local
energy density. In fact, in this setting, the local energy density
at each point is
\begin{equation}
   \mu = \frac{1}{16\pi} R,
\end{equation}
where $R$ is the scalar curvature of the 3-manifold (which has zero second
fundamental form in the spacetime) at each point.  Thus, we
note that $(\real^3,\delta_{ij})$ has zero energy density at each
point as well as zero total mass.
This is appropriate since $(\real^3,\delta_{ij})$ is in
fact a ``$t=0$'' slice of Minkowski spacetime, which represents a vacuum.
Classically, physicists consider $\mu \ge 0$ to be a physical assumption.
Hence, from this point on, we will not only assume that $(M^3,g)$ is
asymptotically flat, but also that it has nonnegative scalar curvature,
\begin{equation}
   R \ge 0.
\end{equation}

This notion of energy density also helps us understand total mass better.
After all, we can take any asymptotically flat manifold and then change the
metric to be perfectly flat outside a large compact set, thereby giving
the new metric zero total mass.  However, if we introduce the physical
condition
that both metrics have nonnegative scalar curvature, then it is a beautiful
theorem that this is in fact not possible, unless the original metric
was already $(\real^3,\delta_{ij})$!  (This theorem is actually a
corollary to the positive mass theorem discussed in a moment.)
Thus, the curvature obstruction of having nonnegative scalar curvature
at each point is a very interesting condition.

Also, notice the indirect connection between the total mass and local
energy density.  At this point, there does not seem to be much of a
connection at all.  Total mass is the $1/r$ rate at which the metric
becomes flat at infinity, and local energy density is the scalar curvature
at each point.  Furthermore, if a metric is changed in a compact set,
local energy density is changed, but the total mass is unaffected.

The reason for this is that the total mass is {\it not} the integral
of the local energy density over the manifold.  In fact, this integral
fails to take potential energy into account (which would be expected to
contribute a negative energy) as well as gravitational energy (discussed
in a moment).  Hence, it is not initially clear what we should expect the
relationship between total mass and local energy density to be, so let us
begin with an example.

\subsection{Example using superharmonic functions in $\real^3$}\vskip -5mm \hspace{5mm}

Once again, let us return to the $(\real^3,u(x)^4\delta_{ij})$ example.
The formula for the scalar curvature is
\begin{equation}
   R = -8 u(x)^{-5} \Delta u(x).
\end{equation}
Hence, since the physical assumption of nonnegative energy density implies
nonnegative scalar curvature, we see that $u(x) > 0$ must be superharmonic
($\Delta u \le 0$).  For simplicity, let's also assume that $u(x)$ is
harmonic outside a bounded set so that we can expand $u(x)$ at infinity using
spherical harmonics.  Hence, $u(x)$ has the asymptotics of
equation \ref{eqn:exp}.
By the maximum principle, it follows that the minimum
value for $u(x)$ must be $a$, referring to equation \ref{eqn:exp}.  Hence,
$b \ge 0$, which implies that $m \ge 0$!  Thus we see that the assumption
of nonnegative energy density at each point of $(\real^3,u(x)^4\delta_{ij})$
implies that the total mass is also nonnegative, which is what one would
hope.

\subsection{The positive mass theorem}\vskip -5mm \hspace{5mm}

More generally, suppose we have any asymptotically flat manifold with
nonnegative scalar curvature, is it true that the total mass is also
nonnegative?  The answer is {\it yes}, and this fact is know as the
positive mass theorem, first proved by Schoen and Yau \cite{SY3} in 1979
using minimal surface techniques and then by Witten \cite{Wi} in 1981 using
spinors.

\begin{theorem} (Schoen-Yau)  Let $(M^3,g)$ be any asymptotically flat,
complete Riemannian
manifold with nonnegative scalar curvature.  Then the total mass
$m \ge 0$, with equality if and only if $(M^3,g)$ is isometric to
$(\real^3,\delta)$.
\end{theorem}

\subsection{Black holes}\vskip -5mm \hspace{5mm}

Another very interesting and natural phenomenon in general relativity is the existence
of black holes.  Instead of thinking of black holes as singularities in a spacetime,
we will think of black holes in terms of their horizons.    Given
a surface in a spacetime, suppose that it admits an outward shell of light.  If the
surface area of this shell of light is decreasing everywhere on the surface, then this
is called a trapped surface.  The outermost boundary of these trapped surfaces is called
the apparent horizon of the black hole.  Apparent horizons can be computed based on
their local geometry, and an apparent horizon always implies the existence of an
event horizon outside of it \cite{HE}.

Now let us return to the case we are considering in this paper where $(M^3,g)$ is a
``t = 0'' slice of a spacetime with zero second fundamental form.  Then it is a very
nice geometric fact that apparent horizons of black holes intersected with
$M^3$ correspond
to the connected components of the outermost minimal surface
$\Sigma_0$ of $(M^3,g)$.

All of the surfaces we are considering
in this paper will be required to be smooth boundaries of open bounded regions, so that
outermost is well-defined with respect to a chosen end of the manifold
\cite{Bray}.  A minimal surface in $(M^3,g)$ is a surface which is a critical point
of the area function with respect to any smooth variation of the surface.  The first
variational calculation implies that minimal surfaces have zero mean curvature.  The
surface $\Sigma_0$ of $(M^3,g)$ is defined as the boundary of the union of the open
regions bounded by all of the minimal surfaces in $(M^3,g)$.  It turns out that
$\Sigma_0$ also has to be a minimal surface, so we call $\Sigma_0$ the {\bf outermost
minimal surface}.

We will also define a surface to be {\bf (strictly) outer minimizing} if every surface
which encloses it has (strictly) greater area.  Note that outermost minimal surfaces
are strictly outer minimizing.  Also, we define a {\bf horizon} in our context to be
any minimal surface which is the boundary of a bounded open region.

It also follows from a
stability argument (using the Gauss-Bonnet theorem interestingly)
that each component of a stable minimal surface (in a 3-manifold with nonnegative
scalar curvature)
must have the topology of a sphere.
Furthermore, there is
a physical argument, based on \cite{P},
which suggests that the mass contributed by the black holes (thought of as the
connected components of $\Sigma_0$) should be defined to be
$\sqrt{A_0/16\pi}$, where $A_0$ is the area of $\Sigma_0$.  Hence, the
physical argument that the total mass should be greater than
or equal to the mass contributed by the black holes yields that following
geometric statement.



\vspace{.1in}\noindent {\bf The Riemannian Penrose Inequality}
\newline {\it Let $(M^3,g)$ be a complete, smooth, 3-manifold with
nonnegative scalar curvature which is harmonically flat at
infinity with total mass $m$ and which has an outermost minimal
surface $\Sigma_0$ of area $A_0$. Then
\begin{equation}\label{eqn:RPI}
m \ge \sqrt{\frac{A_0}{16\pi}},
\end{equation}
with equality if and only if $(M^3,g)$ is isometric to the Schwarzschild metric $(\real^3 \backslash \{0\}$,
$(1+\frac{m}{2|x|})^4\delta_{ij})$ outside their respective outermost minimal surfaces. }\vspace{.1in}

The above statement has been proved by the author \cite{Bray}, and by
Huisken and Ilmanen \cite{HI}
where $A_0$ is defined instead to be the area of the
largest connected component of $\Sigma_0$.  We will discuss both approaches
in this paper, which are very different, although they both involve flowing
surfaces and/or metrics.

We also clarify that the above statement is with respect to a chosen end of
$(M^3,g)$, since both the total mass and the definition of outermost refer to
a particular end.  In fact, nothing very important is gained by considering
manifolds with more than one end, since extra ends can always be compactified
by connect summing them (around a neighborhood of infinity) with large spheres
while still preserving nonnegative scalar curvature, for example.  Hence, we
will typically consider manifolds with just one end.
In the case that the manifold has multiple ends,
we will require every surface (which could have multiple connected
components) in this paper
to enclose all of the ends of the manifold except the chosen end.

Other contributions on the Penrose Conjecture have also been made by
Herzlich \cite{He} using the Dirac
operator which Witten \cite{Wi} used to prove the positive mass theorem,
by Gibbons \cite{Gi} in the special case of collapsing shells,
by Tod \cite{T}, by Bartnik \cite{Ba3} for quasi-spherical metrics, and by
the author \cite{Bray_thesis} using isoperimetric surfaces.  There is also
some interesting work of Ludvigsen and Vickers \cite{LV} using
spinors and Bergqvist \cite{Berg}, both concerning the Penrose inequality for
null slices of a space-time.

\subsection{The Schwarzschild metric}\vskip -5mm \hspace{5mm}

The Schwarzschild metric
$(\real^3 \backslash \{0\}, (1+\frac{m}{2|x|})^4\delta_{ij})$,
referred to in the above
statement of the Riemannian Penrose Inequality, is a particularly important
example to consider, and corresponds to a zero-second fundamental form,
space-like slice of the usual (3+1)-dimensional Schwarzschild metric
(which represents a spherically symmetric static black hole in vacuum).  The
3-dimensional Schwarzschild metrics have total mass $m > 0$
and are characterized
by being the only spherically symmetric,
geodesically complete, zero scalar curvature
3-metrics, other than $(\real^3,\delta_{ij})$.  They can also be embedded in
4-dimensional Euclidean space $(x,y,z,w)$ as the set of points satisfying
$|(x,y,z)| = \frac{w^2}{8m} + 2m$,
which is a parabola rotated around an $S^2$.  This
last picture allows us to see that the Schwarzschild metric, which
has two ends, has a $Z_2$
symmetry which fixes the sphere with $w=0$ and $|(x,y,z)| = 2m$,
which is clearly minimal.  Furthermore, the area of this sphere is
$4\pi (2m)^2$, giving equality in the Riemannian Penrose Inequality.

\section{The conformal flow of metrics}\label{sec:cfm}

\vskip -5mm \hspace{5mm}

Given any initial Riemannian manifold $(M^3,g_0)$
which has nonnegative scalar curvature and
which is harmonically flat at infinity,
we will define a continuous,
one parameter family of metrics $(M^3,g_t)$,
$0 \le t < \infty$.
This family of metrics will converge to a 3-dimensional Schwarzschild
metric and will have other special properties which will
allow us to prove the Riemannian Penrose Inequality for
the original metric $(M^3,g_0)$.

In particular, let $\Sigma_0$ be the outermost minimal surface of $(M^3,g_0)$
with area $A_0$.  Then we will also define a family of surfaces $\Sigma(t)$ with
$\Sigma(0) = \Sigma_0$
such that $\Sigma(t)$ is minimal in $(M^3,g_t)$.  This is natural since
as the metric $g_t$ changes, we expect that the location
of the horizon $\Sigma(t)$ will also change.
Then the interesting quantities to keep track of in this flow are $A(t)$,
the total area of the horizon $\Sigma(t)$ in $(M^3,g_t)$, and
$m(t)$, the total mass of $(M^3,g_t)$ in
the chosen end.

In addition to all of the metrics $g_t$ having nonnegative scalar curvature,
we will also have the very nice properties that
\begin{eqnarray}
A'(t) & = &   0 ,  \\
m'(t) & \le & 0
\end{eqnarray}
for all $t \ge 0$.  Then since $(M^3,g_t)$
converges to a Schwarzschild metric (in an appropriate sense)
which gives equality in the Riemannian Penrose Inequality as described in the
introduction,
\begin{equation}\label{eqn:argument}
m(0) \ge m(\infty) = \sqrt{\frac{A(\infty)}{16\pi}}
     = \sqrt{\frac{A(0)}{16\pi}}
\end{equation}
which proves the Riemannian Penrose Inequality for the original metric
$(M^3,g_0)$.  The hard part, then, is to find a flow of metrics
which preserves nonnegative scalar curvature and the area of the horizon,
decreases total mass, and converges to a Schwarzschild metric as $t$ goes
to infinity.

\subsection{The definition of the flow}\label{sec:def}\vskip -5mm \hspace{5mm}

In fact, the metrics $g_t$ will all be conformal to $g_0$.  This conformal
flow of metrics can be thought of as the solution to a first order o.d.e.~in $t$
defined by equations \ref{eqn:ODE1}, \ref{eqn:ODE2},
\ref{eqn:ODE3}, and \ref{eqn:ODE4}.
Let
\begin{equation}\label{eqn:ODE1}
   g_t = u_t(x)^4 g_0
\end{equation}
and $u_0(x) \equiv 1$.  Given the metric $g_t$, define
\begin{equation}\label{eqn:ODE2}
   \Sigma(t) = \mbox{the outermost minimal area enclosure of }
               \Sigma_0 \mbox{ in } (M^3,g_t)
\end{equation}
where $\Sigma_0$ is the original outer minimizing horizon in $(M^3,g_0)$.
In the cases in which we are interested,
$\Sigma(t)$ will not touch $\Sigma_0$, from which it follows
that $\Sigma(t)$ is actually a strictly outer minimizing horizon of $(M^3,g_t)$.
Then given the horizon $\Sigma(t)$, define $v_t(x)$ such that
\begin{equation}\label{eqn:ODE3}
\left\{
\begin{array}{r l l l}
\Delta_{g_0} v_t(x) & \equiv & 0 & \mbox{ outside } \Sigma(t) \\
v_t(x) & = & 0 & \mbox{ on } \Sigma(t) \\
\lim_{x \rightarrow \infty} v_t(x) & = & -e^{-t} & \\
\end{array}
\right.
\end{equation}
and $v_t(x) \equiv 0$ inside $\Sigma(t)$.  Finally, given $v_t(x)$, define
\begin{equation}\label{eqn:ODE4}
u_t(x) = 1 + \int_0^t v_s(x) ds
\end{equation}
so that $u_t(x)$ is continuous in $t$ and has $u_0(x) \equiv 1$.

Note that equation \ref{eqn:ODE4} implies that the first order rate of change of
$u_t(x)$ is given by $v_t(x)$.  Hence, the first order rate of change of $g_t$
is a function of itself, $g_0$, and $v_t(x)$ which is a function of $g_0$,
$t$, and $\Sigma(t)$ which is in
turn a function of $g_t$ and $\Sigma_0$.  Thus, the first order rate of change
of $g_t$ is a function of $t$, $g_t$, $g_0$, and $\Sigma_0$.

\begin{theorem}\label{thm:existence}
Taken together, equations \ref{eqn:ODE1}, \ref{eqn:ODE2}, \ref{eqn:ODE3},
and \ref{eqn:ODE4} define a first order o.d.e.~in $t$ for $u_t(x)$ which has
a solution which is Lipschitz in the $t$ variable,
$C^1$ in the $x$ variable everywhere,
and smooth in the $x$ variable outside $\Sigma(t)$.
Furthermore, $\Sigma(t)$ is a smooth,
strictly outer minimizing horizon in $(M^3,g_t)$
for all $t \ge 0$, and
$\Sigma(t_2)$ encloses but does not touch $\Sigma(t_1)$
for all $t_2 > t_1 \ge 0$.
\end{theorem}

Since $v_t(x)$ is a superharmonic function in $(M^3,g_0)$
(harmonic everywhere except on $\Sigma(t)$, where it is weakly
superharmonic), it follows that
$u_t(x)$ is superharmonic as well.  Thus, from equation \ref{eqn:ODE4} we see
that $\lim_{x \rightarrow \infty} u_t(x) = e^{-t}$ and consequently that
$u_t(x) > 0$ for all $t$ by the maximum principle.  Then since
\begin{equation}\label{eqn:sc3}
R(g_t) = u_t(x)^{-5}(-8 \Delta_{g_0} + R(g_0))u_t(x),
\end{equation}
it follows that $(M^3,g_t)$ is an asymptotically flat manifold with
nonnegative scalar curvature.

Even so, it still may not seem like $g_t$ is particularly naturally
defined since the rate of change of $g_t$
appears to depend on $t$ and the original metric $g_0$
in equation \ref{eqn:ODE3}.  We would prefer a flow where the rate
of change of $g_t$ can be defined purely as a function of $g_t$ (and $\Sigma_0$
perhaps), and interestingly enough this
actually does turn out to be the case.
In section \ref{sec:harmonic} we prove this very important fact and define
a new equivalence class of metrics called the harmonic conformal class.
Then once we decide to find a flow of metrics which stays inside the
harmonic conformal class of the original metric (outside the horizon)
and keeps the area of the horizon $\Sigma(t)$ constant, then we are
basically forced
to choose the particular conformal flow of metrics defined above.

\begin{theorem}\label{thm:monotone}
The function
$A(t)$ is constant in $t$ and $m(t)$ is non-increasing in $t$,
for all $t \ge 0$.
\end{theorem}

The fact that $A'(t) = 0$ follows from the fact that to first order the
metric is not changing on $\Sigma(t)$ (since $v_t(x) = 0$ there) and from
the fact that to first order the area of $\Sigma(t)$ does not change as it
moves outward since $\Sigma(t)$ is a critical point for area in $(M^3,g_t)$.
Hence, the interesting part of theorem \ref{thm:monotone} is proving that
$m'(t) \le 0$. Curiously, this follows from a nice trick
using the Riemannian positive mass theorem, which we describe in section
\ref{sec:monotonicity}.

Another important aspect of this conformal flow of the metric is that
outside the horizon $\Sigma(t)$, the manifold $(M^3,g_t)$ becomes more
and more spherically symmetric and ``approaches'' a Schwarzschild manifold
$(\real^3 \backslash \{0\}, s)$
in the limit as $t$ goes to
$\infty$.  More precisely,

\begin{theorem}\label{thm:limit}
For sufficiently large $t$, there exists a diffeomorphism $\phi_t$
between \\$(M^3,g_t)$ outside the horizon $\Sigma(t)$ and a fixed
Schwarzschild manifold $(\real^3 \backslash \{0\}, s)$ outside its
horizon.  Furthermore, for all $\epsilon > 0$, there exists a $T$
such that for all $t>T$, the metrics $g_t$ and $\phi^{*}_t(s)$
(when determining the lengths of unit vectors of $(M^3,g_t)$) are
within $\epsilon$ of each other and the total masses of the two
manifolds are within $\epsilon$ of each other.  Hence,
\begin{equation}
\lim_{t \rightarrow \infty} \frac{m(t)}{\sqrt{A(t)}} = \sqrt{\frac{1}{16\pi}}.
\end{equation}
\end{theorem}

Theorem \ref{thm:limit} is not that surprising really although a careful
proof is reasonably long.  However, if one is willing to believe that
the flow of metrics converges to a spherically symmetric metric outside
the horizon, then theorem \ref{thm:limit} follows from two facts.  The
first fact is that the
scalar curvature of $(M^3,g_t)$
eventually becomes identically zero outside the
horizon $\Sigma(t)$ (assuming $(M^3,g_0)$ is harmonically flat).
This follows from the facts that $\Sigma(t)$
encloses any compact set in a finite amount of time, that harmonically
flat manifolds have zero scalar curvature outside a compact set, that
$u_t(x)$ is harmonic outside $\Sigma(t)$, and equation
\ref{eqn:sc3}.  The second fact is that
the Schwarzschild metrics are the only complete,
spherically symmetric 3-manifolds
with zero scalar curvature (except for the flat metric on $R^3$).

The Riemannian Penrose inequality, inequality \ref{eqn:RPI}, then
follows from equation \ref{eqn:argument} using theorems
\ref{thm:existence}, \ref{thm:monotone} and \ref{thm:limit},
for harmonically flat manifolds \cite{Bray}.
Since asymptotically flat manifolds can be approximated arbitrarily well by
harmonically flat manifolds while changing the relevant quantities arbitrarily
little, the asymptotically flat case also follows.
Finally, the case of equality of the Penrose
inequality follows from
a more careful analysis of these same arguments.

\subsection{Qualitative discussion}

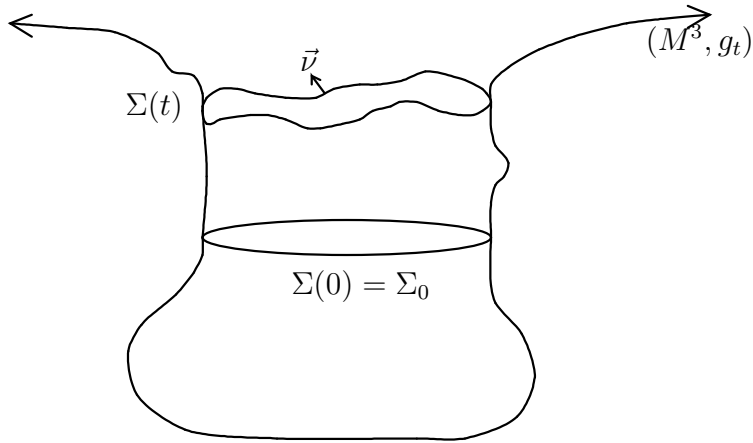
\vspace{.5in}
\begin{center}
\input{040401.eepic}
\end{center}
\vspace{.5in}

The diagrams above and below are meant to help illustrate some of the
properties of the conformal flow of the metric.  The above picture is
the original metric which has a strictly outer minimizing horizon
$\Sigma_0$.  As $t$ increases, $\Sigma(t)$ moves outwards, but never
inwards.  In the diagram below, we can observe one of the consequences
of the fact that $A(t) = A_0$ is constant in $t$.  Since the metric
is not changing inside $\Sigma(t)$, all of the horizons $\Sigma(s)$,
$0 \le s \le t$ have area $A_0$ in $(M^3,g_t)$.  Hence, inside $\Sigma(t)$,
the manifold $(M^3,g_t)$ becomes
cylinder-like
in the sense that it is laminated (meaning foliated but with some gaps allowed)
by all of the previous horizons
which all have the same area $A_0$ with respect to the metric $g_t$.

\newpage

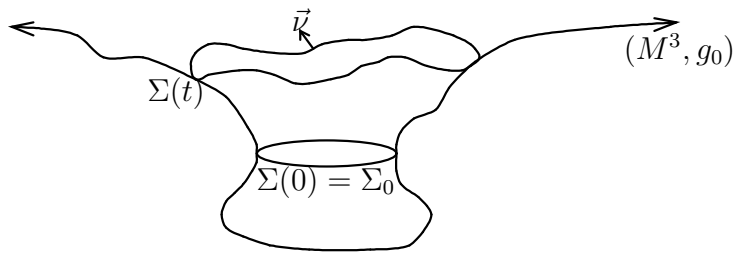
\vspace*{0.5mm}

\begin{center}
\input{040402.eepic}
\end{center}

\vspace*{0.5mm}

Now let us suppose that the
original horizon $\Sigma_0$ of $(M^3,g)$ had two components, for
example.
Then each of the components of the horizon will move outwards as $t$
increases, and at some point before they touch they
will suddenly jump outwards
to form a horizon with a single component
enclosing the previous horizon with two components.  Even horizons with
only one component will sometimes jump outwards, but no more than a countable
number of times.  It is interesting that
this phenomenon of surfaces jumping is also found in the Huisken-Ilmanen
approach to the Penrose conjecture using their generalized $1/H$ flow.

\subsection{Proof that $m'(t) \le 0$}\label{sec:monotonicity}\vskip -5mm \hspace{5mm}

The most surprising aspect of the flow defined in section \ref{sec:def}
is that $m'(t) \le 0$.  As mentioned in that section, this important
fact follows from a nice trick using the Riemannian positive mass theorem.

The first step is to realize that while the rate of change of $g_t$ appears
to depend on $t$ and $g_0$, this is in fact an illusion.  As is described
in detail in section \ref{sec:harmonic}, the rate of change of $g_t$
can be described purely in terms of $g_t$ (and $\Sigma_0$).  It is also
true that the rate of change of $g_t$ depends only on $g_t$ and
$\Sigma(t)$.  Hence, there is no special value of $t$, so proving
$m'(t) \le 0$ is equivalent to proving $m'(0) \le 0$.  Thus, without
loss of generality, we take $t = 0$ for convenience.

Now expand the harmonic function
$v_0(x)$, defined in equation \ref{eqn:ODE3}, using spherical harmonics
at infinity, to get
\begin{equation}\label{eqn:expansion}
   v_0(x) = -1 + \frac{c}{|x|} + {\cal O}\left(\frac{1}{|x|^2}\right)
\end{equation}
for some constant $c$.  Since the rate of change of the metric $g_t$ at
$t=0$ is given by $v_0(x)$ and since the total mass $m(t)$ depends on the
$1/r$ rate at which the metric $g_t$ becomes flat at infinity (see
equation \ref{eqn:tmass}), it is not surprising that direct calculation
gives us that
\begin{equation}
   m'(0) = 2(c - m(0)).
\end{equation}
Hence, to show that $m'(0) \le 0$, we need to show that
\begin{equation}\label{eqn:conj}
c \le m(0).
\end{equation}

In fact, counterexamples to equation \ref{eqn:conj} can be found if we
remove either of the requirements that $\Sigma(0)$ (which is used in the
definition of $v_0(x)$) be a minimal surface or that $(M^3,g_0)$ have
nonnegative scalar curvature.  Hence, we quickly see that equation
\ref{eqn:conj} is a fairly deep conjecture which says something quite
interesting about manifold with nonnegative scalar curvature.  Well,
the Riemannian positive mass theorem is also a deep conjecture which says
something quite interesting about manifolds with nonnegative scalar
curvature.  Hence, it is natural to try to use the Riemannian positive
mass theorem to prove equation \ref{eqn:conj}.

Thus, we want to create a manifold whose total mass depends on $c$
from equation \ref{eqn:expansion}.
The idea is to use a reflection trick similar to one used by Bunting and
Masood-ul-Alam for another purpose in \cite{BM}.  First, remove the region
of $M^3$ inside $\Sigma(0)$ and then reflect the remainder of $(M^3,g_0)$
through $\Sigma(0)$.  Define the resulting Riemannian manifold to be
$(\bar{M}^3, \bar{g}_0)$ which has two asymptotically flat ends since
$(M^3, g_0)$ has exactly one asymptotically flat end not contained by
$\Sigma(0)$.  Note that $(\bar{M}^3, \bar{g}_0)$ has nonnegative scalar
curvature everywhere except on $\Sigma(0)$ where the metric has corners.
In fact, the fact that $\Sigma(0)$ has zero mean curvature (since it is a
minimal surface) implies that $(\bar{M}^3, \bar{g}_0)$ has
{\it distributional} nonnegative scalar curvature everywhere, even on
$\Sigma(0)$.  This notion is made rigorous in \cite{Bray}.  Thus we have
used the fact that $\Sigma(0)$ is minimal in a critical way.

Recall from equation \ref{eqn:ODE3} that $v_0(x)$ was defined to be the
harmonic function equal to zero on $\Sigma(0)$ which goes to $-1$
at infinity.  We want to reflect $v_0(x)$ to be defined on all of
$(\bar{M}^3, \bar{g}_0)$.  The trick here is to define $v_0(x)$ on
$(\bar{M}^3, \bar{g}_0)$ to be the harmonic function which goes to
$-1$ at infinity in the original end and goes to $1$ at infinity in
the reflect end.  By symmetry, $v_0(x)$ equals $0$ on $\Sigma(0)$ and so
agrees with its original definition on $(M^3,g_0)$.

The next step is to compactify one end of $(\bar{M}^3, \bar{g}_0)$.
By the maximum principle, we know that $v_0(x) > -1$ and $c>0$, so the
new Riemannian manifold $(\bar{M}^3, (v_0(x) + 1)^4 \bar{g}_0)$ does the
job quite nicely and compactifies the original end to a point.  In fact,
the compactified point at infinity and the metric there
can be filled in smoothly
(using the fact that $(M^3,g_0)$ is harmonically flat).  It then follows
from equation \ref{eqn:sc3} that this new compactified manifold has
nonnegative scalar curvature since $v_0(x) + 1$ is harmonic.

The last step is simply to apply the Riemannian positive mass theorem
to $(\bar{M}^3, (v_0(x) + 1)^4 \bar{g}_0)$.  It is not surprising that
the total mass $\tilde{m}(0)$ of this manifold involves $c$, but it is quite
lucky that direct calculation yields
\begin{equation}
   \tilde{m}(0) = -4 (c - m(0)),
\end{equation}
which must be positive by the Riemannian positive mass theorem.  Thus,
we have that
\begin{equation}
   m'(0) = 2 (c - m(0)) = -\frac12 \tilde{m}(0) \le 0.
\end{equation}



\subsection{The harmonic conformal class of a metric} \label{sec:harmonic}\vskip -5mm \hspace{5mm}

As a final topic which is also of independent interest,
we define a new equivalence class
and partial ordering of conformal metrics.  These new objects
provide a natural motivation for studying
conformal flows of metrics to try to prove the Riemannian Penrose inequality.
Let
\begin{equation}\label{eqn:conf}
g_2 = u(x)^{\frac{4}{n-2}} g_1,
\end{equation}
where $g_2$ and $g_1$ are metrics on
an $n$-dimensional manifold $M^n$, $n \ge 3$.  Then we get the surprisingly
simple identity that
\begin{equation}\label{eqn:identity}
\Delta_{g_1}(u \phi) = u^\frac{n+2}{n-2} \Delta_{g_2}(\phi)
                       + \phi \Delta_{g_1}(u)
\end{equation}
for any smooth function $\phi$.
This motivates us to define the following relation.
\begin{definition}
Define
\[ g_2 \sim g_1 \]
if and only if
equation \ref{eqn:conf} is satisfied with $\Delta_{g_1}(u) = 0$ and $u(x) > 0$.
\end{definition}
Then from equation \ref{eqn:identity} we get the following lemma.
\begin{lemma}
The relation $\sim$ is reflexive, symmetric, and transitive, and hence
is an equivalence relation.
\end{lemma}
Thus, we can define the following equivalence class of metrics.
\begin{definition}
Define
\[ [g]_H = \{\bar{g} \;\;|\;\; \bar{g} \sim g \} \]
to be the {\bf harmonic conformal class} of the metric $g$.
\end{definition}
Of course, this definition is most interesting when $(M^n,g)$ has
nonconstant positive harmonic functions, which happens for example
when $(M^n,g$) has a boundary.

Also, we can modify the relation $\sim$ to get another relation $\succeq$.
\begin{definition}
Define
\[ g_2 \succeq g_1 \]
if and only if equation \ref{eqn:conf} is satisfied
with $- \Delta_{g_1}(u) \ge 0$ and $u(x) > 0$.
\end{definition}
Then from equation \ref{eqn:identity} we get the following lemma.
\begin{lemma}
The relation $\succeq$ is reflexive and transitive, and hence
is a partial ordering.
\end{lemma}
Since $\succeq$ is defined in terms of superharmonic functions, we
will call it the superharmonic partial ordering of metrics on $M^n$.
Then it is natural to define the following set of metrics.
\begin{definition}
Define
\[ [g]_S = \{\bar{g} \;\;|\;\; \bar{g} \succeq g \} .\]
\end{definition}
This set of metrics has the property that if $\bar{g} \in [g]_S$, then
$[\bar{g}]_S \subset [g]_S$

Also, the scalar curvature transforms nicely under a conformal change of
the metric.  In fact, assuming equation \ref{eqn:conf} again,
\begin{equation}\label{eqn:scalar_curv}
  R(g_2) = u(x)^{-(\frac{n+2}{n-2})}\left(-c_n \Delta_{g_1} + R(g_1)\right)u(x)
\end{equation}
where $c_n = \frac{4(n-1)}{n-2}$.  This gives us the following lemma.
\begin{lemma}
The sign of the scalar curvature is preserved pointwise by $\sim$.  That is,
if $g_2 \sim g_1$, then $sgn(R(g_2)(x)) = sgn(R(g_1)(x))$ for all $x \in M^n$.
Also, if $g_2 \succeq g_1$, and $g_1$ has non-negative scalar curvature, then
$g_2$ has non-negative scalar curvature.
\end{lemma}

Hence, the harmonic conformal equivalence relation $\sim$ and the superharmonic
partial ordering $\succeq$ are useful for studying questions about
scalar curvature.  In particular, these notions are useful for studying the
Riemannian Penrose inequality which concerns asymptotically flat 3-manifolds
$(M^3,g)$ with non-negative scalar curvature.  Given such a manifold, define
$m(g)$ to be the total mass of $(M^3,g)$ and $A(g)$ to be the area of the
outermost horizon (which could have multiple components) of $(M^3,g)$.  Define
$P(g) = \frac{m(g)}{\sqrt{A(g)}}$ to be the Penrose quotient of $(M^3,g)$.
Then an interesting question is to ask which metric in $[g]_S$ minimizes
$P(g)$.

Section \ref{sec:cfm} of this
paper can be viewed as an answer to the above question.
We showed that
there exists a conformal flow of metrics (starting with $g_0$)
for which the Penrose quotient
was non-increasing, and in fact this conformal flow stays inside $[g_0]_S$.
Furthermore, $g_{t_2} \in [g_{t_1}]_S$ for all $t_2 \ge t_1 \ge 0$.  We showed
that no matter which metric we start with,
the metric converges to a Schwarzschild metric
outside its horizon.  Hence, the minimum value of $P(g)$ in $[g]_S$ is
achieved in the limit
by metrics converging to a Schwarzschild metric (outside their
respective horizons).

In the case that $g$ is harmonically flat at infinity, a Schwarzschild
metric (outside the horizon)
is contained in $[g]_S$.  More generally, given any asymptotically
flat manifold $(M^3,g)$, we can use $\real^3 \backslash  B_r(0)$ as a
coordinate chart for the asymptotically flat end of $(M^3,g)$ which we are
interested in, where the metric $g_{ij}$ approaches $\delta_{ij}$
at infinity in this
coordinate chart.
Then we can consider the conformal metric
\begin{equation}
   g_C = \left(1 + \frac{C}{|x|}\right)^4 g
\end{equation}
in this end.  In the limit as $C$ goes to infinity,
the horizon will approach the coordinate sphere of radius $C$.  Then
outside this horizon in the limit as $C$ goes to infinity, the function
$(1 + \frac{C}{|x|})$ will be close to a superharmonic function
on $(M^3,g)$ and the
metric $g_C$ will approach a Schwarzschild metric
(since the metric $g$ is approaching the standard metric on
$\real^3$).  Hence,
the Penrose quotient of $g_C$ will approach $(16 \pi)^{-1/2}$, which is the
Penrose quotient of a Schwarzschild metric.

As a final note, we prove that the first order o.d.e.~for $\{g_t\}$ defined in equations \ref{eqn:ODE1},
\ref{eqn:ODE2}, \ref{eqn:ODE3}, and \ref{eqn:ODE4} is naturally defined in the sense that the rate of change of
$g_t$ is a function only of $g_t$ and not of $g_0$ or $t$. To see this, given any solution $g_t = u_t(x)^4 g_0$ to
equations \ref{eqn:ODE1}, \ref{eqn:ODE2}, \ref{eqn:ODE3}, and \ref{eqn:ODE4}, choose any $s > 0$ and define
$\bar{u}_t(x) = u_t(x)/u_s(x)$ so that
\begin{equation}\label{eqn:ODE1_bar}
   g_t = \bar{u}_t(x)^4 g_s
\end{equation}
and $\bar{u}_s(x) \equiv 1$.
Then  define $\bar{v}_t(x)$ such that
\begin{equation}\label{eqn:ODE3_bar}
\left\{
\begin{array}{r l l l}
\Delta_{g_s} \bar{v}_t(x) & \equiv & 0 & \mbox{ outside } \Sigma(t) \\
\bar{v}_t(x) & = & 0 & \mbox{ on } \Sigma(t) \\
\lim_{x \rightarrow \infty} \bar{v}_t(x) & = & -e^{-(t-s)} & \\
\end{array}
\right.
\end{equation}
and $\bar{v}_t(x) \equiv 0$ inside $\Sigma(t)$.
Then what we want to show is
\begin{equation}\label{eqn:ODE4_bar}
\bar{u}_t(x) = 1 + \int_s^t \bar{v}_r(x) dr
\end{equation}
To prove the above equation, we observe that
from equations \ref{eqn:identity}, \ref{eqn:ODE3_bar}, and \ref{eqn:ODE3}
it follows that
\begin{equation}
   v_t(x) = \bar{v}_t(x) \, u_s(x)
\end{equation}
since $\lim_{x \rightarrow \infty} u_s(x) = e^{-s}$.
Hence, since
\begin{equation}
u_t(x) = u_s(x) + \int_s^t v_r(x) dr
\end{equation}
by equation \ref{eqn:ODE4}, dividing through by $u_s(x)$ yields
equation \ref{eqn:ODE4_bar} as desired.
Thus, we see that
the rate of change of $g_t(x)$ at
$t=s$ is a function of $\bar{v}_s(x)$ which in turn
is just a function of $g_s(x)$ and the horizon $\Sigma(s)$.  Hence, to
understand properties of the flow we need only analyze the behavior of the
flow for $t$ close to zero, since any metric in the flow may be chosen to be
the base metric.


\label{lastpage}

\end{document}

%% file: 040401.eepic
\setlength{\unitlength}{0.0003in}
\begingroup\makeatletter\ifx\SetFigFont\undefined
\def\x#1#2#3#4#5#6#7\relax{\def\x{#1#2#3#4#5#6}}%
\expandafter\x\fmtname xxxxxx\relax \def\y{splain}%
\ifx\x\y   
\gdef\SetFigFont#1#2#3{%
  \ifnum #1<17\tiny\else \ifnum #1<20\small\else
  \ifnum #1<24\normalsize\else \ifnum #1<29\large\else
  \ifnum #1<34\Large\else \ifnum #1<41\LARGE\else
     \huge\fi\fi\fi\fi\fi\fi
  \csname #3\endcsname}%
\else
\gdef\SetFigFont#1#2#3{\begingroup
  \count@#1\relax \ifnum 25<\count@\count@25\fi
  \def\x{\endgroup\@setsize\SetFigFont{#2pt}}%
  \expandafter\x
    \csname \romannumeral\the\count@ pt\expandafter\endcsname
    \csname @\romannumeral\the\count@ pt\endcsname
  \csname #3\endcsname}%
\fi
\fi\endgroup
\begin{picture}(12142,7479)(0,-10)
\thicklines
\put(5787,3477){\ellipse{4950}{600}}
\path(5412,5952)(5187,6252)(5187,6102)
\path(5187,6252)(5337,6177)
\path(387,7302)(12,7152)(387,6927)
\path(11637,7452)(12012,7302)(11637,7077)
\path(3312,5727)	(3323.816,5762.782)
	(3342.754,5805.233)
	(3387.000,5877.000)

\path(3387,5877)	(3431.598,5918.727)
	(3491.625,5963.333)
	(3555.590,6002.272)
	(3612.000,6027.000)

\path(3612,6027)	(3679.619,6037.418)
	(3720.437,6038.848)
	(3763.174,6038.385)
	(3805.764,6036.499)
	(3846.142,6033.660)
	(3912.000,6027.000)

\path(3912,6027)	(3979.152,6013.495)
	(4019.917,6003.001)
	(4062.686,5991.307)
	(4105.370,5979.388)
	(4145.878,5968.216)
	(4212.000,5952.000)

\path(4212,5952)	(4263.735,5942.441)
	(4327.160,5931.664)
	(4398.414,5920.303)
	(4435.771,5914.603)
	(4473.637,5908.995)
	(4511.531,5903.559)
	(4548.969,5898.375)
	(4620.547,5889.079)
	(4684.511,5881.742)
	(4737.000,5877.000)

\path(4737,5877)	(4788.935,5873.351)
	(4852.238,5869.288)
	(4923.182,5865.539)
	(4998.038,5862.833)
	(5035.767,5862.098)
	(5073.076,5861.897)
	(5144.568,5863.461)
	(5208.786,5868.252)
	(5262.000,5877.000)

\path(5262,5877)	(5331.783,5904.911)
	(5372.041,5926.278)
	(5413.710,5949.911)
	(5455.165,5973.806)
	(5494.782,5995.956)
	(5562.000,6027.000)

\path(5562,6027)	(5606.191,6039.343)
	(5660.443,6051.426)
	(5721.448,6062.954)
	(5785.897,6073.635)
	(5850.485,6083.174)
	(5911.903,6091.277)
	(5966.844,6097.650)
	(6012.000,6102.000)

\path(6012,6102)	(6078.188,6103.597)
	(6118.993,6102.998)
	(6162.000,6102.000)
	(6205.007,6101.002)
	(6245.812,6100.403)
	(6312.000,6102.000)

\path(6312,6102)	(6357.157,6106.538)
	(6412.098,6113.370)
	(6473.517,6122.039)
	(6538.106,6132.090)
	(6602.557,6143.067)
	(6663.561,6154.515)
	(6717.811,6165.978)
	(6762.000,6177.000)

\path(6762,6177)	(6799.744,6190.325)
	(6844.216,6209.885)
	(6893.356,6233.336)
	(6945.105,6258.337)
	(6997.403,6282.547)
	(7048.192,6303.622)
	(7095.410,6319.220)
	(7137.000,6327.000)

\path(7137,6327)	(7209.731,6326.803)
	(7250.914,6322.991)
	(7294.576,6317.055)
	(7340.168,6309.264)
	(7387.142,6299.881)
	(7434.947,6289.174)
	(7483.035,6277.406)
	(7530.856,6264.845)
	(7577.861,6251.755)
	(7623.501,6238.403)
	(7667.227,6225.054)
	(7708.488,6211.973)
	(7746.737,6199.427)
	(7812.000,6177.000)

\path(7812,6177)	(7880.814,6149.060)
	(7922.284,6130.247)
	(7965.521,6109.515)
	(8008.318,6087.843)
	(8048.468,6066.213)
	(8112.000,6027.000)

\path(8112,6027)	(8151.615,5991.615)
	(8187.000,5952.000)

\path(8187,5952)	(8236.481,5885.602)
	(8257.280,5844.104)
	(8262.000,5802.000)

\path(8262,5802)	(8239.493,5753.365)
	(8197.935,5711.887)
	(8150.909,5677.966)
	(8112.000,5652.000)

\path(8112,5652)	(8039.858,5609.362)
	(7997.641,5589.419)
	(7962.000,5577.000)

\path(7962,5577)	(7911.696,5570.991)
	(7848.938,5570.381)
	(7786.461,5573.081)
	(7737.000,5577.000)

\path(7737,5577)	(7699.161,5581.689)
	(7653.261,5588.858)
	(7602.040,5597.936)
	(7548.236,5608.354)
	(7494.591,5619.539)
	(7443.843,5630.923)
	(7398.733,5641.933)
	(7362.000,5652.000)

\path(7362,5652)	(7324.176,5665.303)
	(7279.081,5684.055)
	(7229.058,5706.309)
	(7176.450,5730.120)
	(7123.598,5753.541)
	(7072.844,5774.625)
	(7026.531,5791.427)
	(6987.000,5802.000)

\path(6987,5802)	(6920.539,5810.042)
	(6879.974,5812.337)
	(6837.293,5813.329)
	(6794.575,5812.905)
	(6753.900,5810.951)
	(6687.000,5802.000)

\path(6687,5802)	(6651.976,5788.232)
	(6610.789,5766.041)
	(6537.000,5727.000)

\path(6537,5727)	(6471.232,5708.973)
	(6430.216,5700.156)
	(6386.884,5691.270)
	(6343.567,5682.162)
	(6302.594,5672.682)
	(6237.000,5652.000)

\path(6237,5652)	(6185.438,5621.303)
	(6124.500,5577.000)
	(6063.562,5532.697)
	(6012.000,5502.000)

\path(6012,5502)	(5946.654,5480.007)
	(5905.898,5469.380)
	(5862.776,5459.235)
	(5819.558,5449.748)
	(5778.510,5441.095)
	(5712.000,5427.000)

\path(5712,5427)	(5675.063,5418.641)
	(5629.857,5408.433)
	(5579.097,5397.154)
	(5525.497,5385.585)
	(5471.773,5374.505)
	(5420.639,5364.695)
	(5374.810,5356.933)
	(5337.000,5352.000)

\path(5337,5352)	(5287.456,5348.421)
	(5224.841,5346.289)
	(5162.055,5347.012)
	(5112.000,5352.000)

\path(5112,5352)	(5072.426,5362.492)
	(5026.007,5379.095)
	(4975.123,5399.935)
	(4922.152,5423.134)
	(4869.476,5446.816)
	(4819.472,5469.105)
	(4774.520,5488.126)
	(4737.000,5502.000)

\path(4737,5502)	(4671.700,5521.680)
	(4631.227,5532.996)
	(4588.417,5544.386)
	(4545.431,5555.165)
	(4504.426,5564.650)
	(4437.000,5577.000)

\path(4437,5577)	(4370.634,5581.271)
	(4329.812,5582.019)
	(4286.820,5582.014)
	(4243.850,5581.381)
	(4203.096,5580.250)
	(4137.000,5577.000)

\path(4137,5577)	(4084.502,5572.452)
	(4020.517,5565.586)
	(3948.913,5556.872)
	(3911.463,5551.969)
	(3873.559,5546.779)
	(3835.684,5541.362)
	(3798.322,5535.775)
	(3727.071,5524.332)
	(3663.674,5512.917)
	(3612.000,5502.000)

\path(3612,5502)	(3566.044,5479.190)
	(3510.525,5444.329)
	(3449.494,5419.553)
	(3387.000,5427.000)

\path(3387,5427)	(3341.731,5471.826)
	(3321.908,5534.347)
	(3315.880,5599.446)
	(3312.000,5652.000)

\path(3312,5652)	(3309.375,5689.384)
	(3312.000,5727.000)

\path(8262,5802)	(8254.381,5861.096)
	(8248.939,5912.263)
	(8245.673,5956.382)
	(8244.585,5994.330)
	(8248.939,6055.232)
	(8262.000,6102.000)

\path(8262,6102)	(8284.243,6142.783)
	(8316.313,6185.045)
	(8355.476,6227.577)
	(8398.999,6269.171)
	(8444.147,6308.619)
	(8488.185,6344.712)
	(8562.000,6402.000)

\path(8562,6402)	(8611.506,6436.990)
	(8673.530,6476.606)
	(8744.116,6518.895)
	(8781.382,6540.430)
	(8819.306,6561.900)
	(8857.392,6583.061)
	(8895.146,6603.668)
	(8967.679,6642.244)
	(9032.949,6675.673)
	(9087.000,6702.000)

\path(9087,6702)	(9145.277,6727.831)
	(9217.255,6757.405)
	(9256.988,6773.113)
	(9298.473,6789.178)
	(9341.152,6805.407)
	(9384.469,6821.606)
	(9427.864,6837.584)
	(9470.781,6853.146)
	(9512.662,6868.101)
	(9552.948,6882.255)
	(9626.508,6907.387)
	(9687.000,6927.000)

\path(9687,6927)	(9728.190,6939.516)
	(9775.245,6953.261)
	(9827.339,6968.017)
	(9883.644,6983.565)
	(9943.336,6999.684)
	(10005.586,7016.156)
	(10069.569,7032.761)
	(10134.458,7049.280)
	(10199.426,7065.493)
	(10263.647,7081.181)
	(10326.294,7096.125)
	(10386.542,7110.105)
	(10443.563,7122.902)
	(10496.530,7134.296)
	(10544.618,7144.069)
	(10587.000,7152.000)

\path(10587,7152)	(10633.890,7159.983)
	(10683.969,7167.894)
	(10737.756,7175.791)
	(10795.775,7183.726)
	(10858.547,7191.755)
	(10926.594,7199.934)
	(11000.438,7208.316)
	(11039.697,7212.601)
	(11080.601,7216.958)
	(11123.215,7221.392)
	(11167.605,7225.913)
	(11213.835,7230.525)
	(11261.970,7235.237)
	(11312.077,7240.054)
	(11364.220,7244.984)
	(11418.465,7250.034)
	(11474.877,7255.210)
	(11533.520,7260.520)
	(11594.461,7265.970)
	(11657.764,7271.567)
	(11723.495,7277.318)
	(11791.719,7283.230)
	(11862.501,7289.310)
	(11935.906,7295.564)
	(11973.613,7298.759)
	(12012.000,7302.000)

\path(3312,5652)	(3316.115,5712.064)
	(3319.054,5763.922)
	(3320.817,5808.456)
	(3321.405,5846.542)
	(3319.054,5906.891)
	(3312.000,5952.000)

\path(3312,5952)	(3291.313,6024.151)
	(3276.267,6066.481)
	(3258.352,6110.025)
	(3237.769,6152.566)
	(3214.716,6191.887)
	(3162.000,6252.000)

\path(3162,6252)	(3124.633,6271.923)
	(3079.040,6283.932)
	(3027.976,6290.477)
	(2974.196,6294.007)
	(2920.454,6296.975)
	(2869.504,6301.829)
	(2824.101,6311.021)
	(2787.000,6327.000)

\path(2787,6327)	(2726.326,6387.952)
	(2697.410,6429.756)
	(2669.273,6474.945)
	(2641.789,6520.390)
	(2614.833,6562.966)
	(2562.000,6627.000)

\path(2562,6627)	(2492.799,6678.290)
	(2451.213,6705.325)
	(2405.865,6732.906)
	(2357.423,6760.744)
	(2306.558,6788.554)
	(2253.940,6816.048)
	(2200.237,6842.940)
	(2146.122,6868.943)
	(2092.262,6893.771)
	(2039.329,6917.137)
	(1987.991,6938.754)
	(1938.920,6958.336)
	(1892.784,6975.595)
	(1850.254,6990.245)
	(1812.000,7002.000)

\path(1812,7002)	(1773.635,7011.535)
	(1729.966,7020.022)
	(1681.754,7027.547)
	(1629.756,7034.197)
	(1574.733,7040.056)
	(1517.444,7045.212)
	(1458.649,7049.751)
	(1399.106,7053.757)
	(1339.576,7057.319)
	(1280.817,7060.521)
	(1223.590,7063.449)
	(1168.653,7066.190)
	(1116.766,7068.829)
	(1068.689,7071.453)
	(1025.180,7074.148)
	(987.000,7077.000)

\path(987,7077)	(920.465,7082.499)
	(883.590,7085.486)
	(843.832,7088.667)
	(800.836,7092.069)
	(754.244,7095.722)
	(703.699,7099.651)
	(648.844,7103.884)
	(589.322,7108.449)
	(524.776,7113.373)
	(454.849,7118.684)
	(417.756,7121.493)
	(379.184,7124.409)
	(339.088,7127.435)
	(297.423,7130.576)
	(254.146,7133.833)
	(209.210,7137.211)
	(162.573,7140.714)
	(114.188,7144.344)
	(64.012,7148.105)
	(12.000,7152.000)

\path(3312,5652)	(3317.206,5592.134)
	(3322.193,5534.378)
	(3326.964,5478.679)
	(3331.523,5424.986)
	(3335.874,5373.248)
	(3340.019,5323.413)
	(3343.963,5275.429)
	(3347.708,5229.246)
	(3351.258,5184.812)
	(3354.617,5142.074)
	(3357.787,5100.982)
	(3360.773,5061.485)
	(3366.205,4987.066)
	(3370.939,4918.406)
	(3375.003,4855.093)
	(3378.425,4796.714)
	(3381.231,4742.858)
	(3383.450,4693.113)
	(3385.109,4647.067)
	(3386.235,4604.307)
	(3386.857,4564.422)
	(3387.000,4527.000)

\path(3387,4527)	(3385.943,4454.828)
	(3384.702,4414.912)
	(3382.950,4371.929)
	(3380.657,4325.496)
	(3377.799,4275.227)
	(3374.346,4220.738)
	(3370.271,4161.645)
	(3365.548,4097.563)
	(3360.147,4028.108)
	(3357.185,3991.245)
	(3354.043,3952.894)
	(3350.718,3913.008)
	(3347.207,3871.538)
	(3343.507,3828.437)
	(3339.613,3783.656)
	(3335.522,3737.146)
	(3331.231,3688.861)
	(3326.737,3638.752)
	(3322.036,3586.771)
	(3317.125,3532.870)
	(3312.000,3477.000)

\path(8262,5802)	(8260.320,5754.984)
	(8258.863,5711.257)
	(8257.631,5670.652)
	(8256.623,5633.005)
	(8255.278,5565.928)
	(8254.830,5508.705)
	(8255.278,5460.019)
	(8256.623,5418.552)
	(8262.000,5352.000)

\path(8262,5352)	(8272.954,5284.469)
	(8281.216,5243.372)
	(8290.890,5200.301)
	(8301.640,5157.443)
	(8313.131,5116.983)
	(8337.000,5052.000)

\path(8337,5052)	(8385.930,4986.158)
	(8457.379,4912.035)
	(8493.279,4872.888)
	(8524.888,4832.895)
	(8548.897,4792.463)
	(8562.000,4752.000)

\path(8562,4752)	(8559.191,4695.020)
	(8539.980,4632.334)
	(8513.029,4573.230)
	(8487.000,4527.000)

\path(8487,4527)	(8455.793,4491.533)
	(8411.569,4452.743)
	(8367.560,4413.580)
	(8337.000,4377.000)

\path(8337,4377)	(8310.587,4312.008)
	(8298.990,4271.548)
	(8288.629,4228.691)
	(8279.626,4185.623)
	(8272.106,4144.529)
	(8262.000,4077.000)

\path(8262,4077)	(8257.816,4035.685)
	(8254.828,3988.262)
	(8253.035,3932.972)
	(8252.438,3868.057)
	(8252.587,3831.442)
	(8253.035,3791.761)
	(8253.782,3748.795)
	(8254.828,3702.324)
	(8256.173,3652.129)
	(8257.816,3597.990)
	(8259.759,3539.687)
	(8262.000,3477.000)

\path(3312,3477)	(3316.115,3416.936)
	(3319.054,3365.077)
	(3320.817,3320.544)
	(3321.405,3282.457)
	(3319.054,3222.109)
	(3312.000,3177.000)

\path(3312,3177)	(3286.546,3105.898)
	(3267.329,3064.486)
	(3245.640,3021.773)
	(3222.871,2979.712)
	(3200.414,2940.261)
	(3162.000,2877.000)

\path(3162,2877)	(3135.086,2837.347)
	(3100.655,2790.229)
	(3060.952,2738.348)
	(3018.225,2684.404)
	(2974.720,2631.096)
	(2932.683,2581.126)
	(2894.361,2537.194)
	(2862.000,2502.000)

\path(2862,2502)	(2804.589,2447.571)
	(2731.101,2385.212)
	(2690.055,2351.798)
	(2647.063,2317.289)
	(2602.813,2281.980)
	(2557.999,2246.167)
	(2513.309,2210.146)
	(2469.435,2174.212)
	(2427.067,2138.662)
	(2386.897,2103.789)
	(2315.910,2037.263)
	(2262.000,1977.000)

\path(2262,1977)	(2232.007,1935.712)
	(2198.360,1884.753)
	(2163.154,1826.821)
	(2128.481,1764.615)
	(2096.436,1700.832)
	(2069.112,1638.170)
	(2048.602,1579.326)
	(2037.000,1527.000)

\path(2037,1527)	(2033.172,1488.276)
	(2031.445,1445.723)
	(2031.746,1399.928)
	(2034.003,1351.478)
	(2038.141,1300.963)
	(2044.089,1248.969)
	(2051.773,1196.083)
	(2061.120,1142.895)
	(2072.057,1089.991)
	(2084.511,1037.960)
	(2098.409,987.388)
	(2113.677,938.864)
	(2130.244,892.976)
	(2148.035,850.310)
	(2187.000,777.000)

\path(2187,777)	(2217.746,733.661)
	(2254.714,689.406)
	(2297.164,644.598)
	(2344.358,599.599)
	(2395.557,554.774)
	(2450.025,510.485)
	(2507.022,467.095)
	(2565.810,424.969)
	(2625.651,384.468)
	(2685.806,345.957)
	(2745.538,309.797)
	(2804.108,276.354)
	(2860.777,245.989)
	(2914.807,219.067)
	(2965.461,195.949)
	(3012.000,177.000)

\path(3012,177)	(3086.020,152.680)
	(3126.825,141.710)
	(3169.922,131.493)
	(3215.135,121.999)
	(3262.286,113.197)
	(3311.198,105.055)
	(3361.692,97.545)
	(3413.592,90.634)
	(3466.720,84.293)
	(3520.899,78.491)
	(3575.950,73.197)
	(3631.697,68.381)
	(3687.963,64.012)
	(3744.568,60.060)
	(3801.338,56.494)
	(3858.092,53.283)
	(3914.655,50.398)
	(3970.848,47.806)
	(4026.495,45.479)
	(4081.418,43.385)
	(4135.438,41.493)
	(4188.379,39.774)
	(4240.064,38.196)
	(4290.314,36.729)
	(4338.953,35.343)
	(4385.803,34.006)
	(4430.685,32.689)
	(4473.424,31.360)
	(4513.841,29.989)
	(4587.000,27.000)

\path(4587,27)	(4640.513,24.712)
	(4698.190,22.714)
	(4759.753,20.992)
	(4824.924,19.532)
	(4893.426,18.320)
	(4964.981,17.342)
	(5039.311,16.582)
	(5077.429,16.280)
	(5116.138,16.027)
	(5155.401,15.822)
	(5195.184,15.663)
	(5235.453,15.547)
	(5276.173,15.475)
	(5317.308,15.442)
	(5358.825,15.448)
	(5400.689,15.492)
	(5442.864,15.570)
	(5485.316,15.681)
	(5528.011,15.825)
	(5570.914,15.998)
	(5613.989,16.198)
	(5657.203,16.426)
	(5700.520,16.677)
	(5743.906,16.951)
	(5787.326,17.246)
	(5830.746,17.560)
	(5874.130,17.892)
	(5917.444,18.238)
	(5960.653,18.599)
	(6003.723,18.971)
	(6046.618,19.354)
	(6089.305,19.745)
	(6131.748,20.142)
	(6173.912,20.544)
	(6215.764,20.949)
	(6257.267,21.355)
	(6298.388,21.761)
	(6339.092,22.164)
	(6379.343,22.564)
	(6419.108,22.957)
	(6458.352,23.342)
	(6497.039,23.718)
	(6535.136,24.083)
	(6609.417,24.771)
	(6680.918,25.393)
	(6749.361,25.933)
	(6814.469,26.378)
	(6875.963,26.714)
	(6933.566,26.926)
	(6987.000,27.000)

\path(6987,27)	(7025.117,26.576)
	(7068.495,25.417)
	(7116.392,23.693)
	(7168.063,21.573)
	(7222.765,19.227)
	(7279.756,16.825)
	(7338.291,14.536)
	(7397.629,12.529)
	(7457.025,10.974)
	(7515.736,10.042)
	(7573.019,9.900)
	(7628.130,10.720)
	(7680.328,12.670)
	(7728.867,15.920)
	(7773.006,20.641)
	(7812.000,27.000)

\path(7812,27)	(7857.638,36.729)
	(7908.503,48.750)
	(7963.779,62.963)
	(8022.648,79.268)
	(8084.293,97.568)
	(8147.899,117.761)
	(8212.646,139.749)
	(8277.720,163.432)
	(8342.302,188.712)
	(8405.576,215.487)
	(8466.725,243.660)
	(8524.932,273.130)
	(8579.380,303.799)
	(8629.252,335.566)
	(8673.731,368.333)
	(8712.000,402.000)

\path(8712,402)	(8764.867,462.620)
	(8790.990,498.877)
	(8816.631,538.401)
	(8841.586,580.667)
	(8865.648,625.154)
	(8888.612,671.338)
	(8910.274,718.695)
	(8930.427,766.702)
	(8948.866,814.836)
	(8965.385,862.574)
	(8979.780,909.392)
	(8991.844,954.767)
	(9001.372,998.175)
	(9008.159,1039.094)
	(9012.000,1077.000)

\path(9012,1077)	(9012.604,1117.363)
	(9009.756,1161.270)
	(9003.776,1208.148)
	(8994.984,1257.422)
	(8983.703,1308.521)
	(8970.252,1360.870)
	(8954.952,1413.896)
	(8938.125,1467.026)
	(8920.091,1519.687)
	(8901.170,1571.304)
	(8881.684,1621.305)
	(8861.953,1669.117)
	(8842.299,1714.166)
	(8823.041,1755.878)
	(8787.000,1827.000)

\path(8787,1827)	(8748.001,1885.213)
	(8693.979,1949.450)
	(8629.982,2017.756)
	(8561.055,2088.180)
	(8492.246,2158.769)
	(8428.603,2227.570)
	(8375.172,2292.632)
	(8337.000,2352.000)

\path(8337,2352)	(8311.992,2416.975)
	(8300.453,2457.434)
	(8289.878,2500.294)
	(8280.504,2543.366)
	(8272.574,2584.466)
	(8262.000,2652.000)

\path(8262,2652)	(8256.248,2708.806)
	(8252.139,2774.011)
	(8250.701,2810.520)
	(8249.674,2850.034)
	(8249.058,2892.857)
	(8248.853,2939.291)
	(8249.058,2989.638)
	(8249.674,3044.199)
	(8250.701,3103.278)
	(8252.139,3167.175)
	(8253.988,3236.194)
	(8256.248,3310.637)
	(8257.532,3349.986)
	(8258.919,3390.805)
	(8260.408,3433.130)
	(8262.000,3477.000)

\put(10737,6702){\makebox(0,0)[lb]{\smash{{{\SetFigFont{12}{14.4}{rm}
$(M^3,g_t)$}}}}}

\put(4812,6327){\makebox(0,0)[lb]{\smash{{{\SetFigFont{12}{14.4}{rm}
$\vec{\nu}$}}}}}

\put(4662,2502){\makebox(0,0)[lb]{\smash{{{\SetFigFont{12}{14.4}{rm}
$\Sigma(0) = \Sigma_0$}}}}}

\put(1812,5577){\makebox(0,0)[lb]{\smash{{{\SetFigFont{12}{14.4}{rm}
$\Sigma(t)$}}}}}

\end{picture}

%% file: 040402.eepic
\setlength{\unitlength}{0.0003in}
\begingroup\makeatletter\ifx\SetFigFont\undefined
\def\x#1#2#3#4#5#6#7\relax{\def\x{#1#2#3#4#5#6}}%
\expandafter\x\fmtname xxxxxx\relax \def\y{splain}%
\ifx\x\y   
\gdef\SetFigFont#1#2#3{%
  \ifnum #1<17\tiny\else \ifnum #1<20\small\else
  \ifnum #1<24\normalsize\else \ifnum #1<29\large\else
  \ifnum #1<34\Large\else \ifnum #1<41\LARGE\else
     \huge\fi\fi\fi\fi\fi\fi
  \csname #3\endcsname}%
\else
\gdef\SetFigFont#1#2#3{\begingroup
  \count@#1\relax \ifnum 25<\count@\count@25\fi
  \def\x{\endgroup\@setsize\SetFigFont{#2pt}}%
  \expandafter\x
    \csname \romannumeral\the\count@ pt\expandafter\endcsname
    \csname @\romannumeral\the\count@ pt\endcsname
  \csname #3\endcsname}%
\fi
\fi\endgroup
\begin{picture}(12081,4020)(0,-10)
\thicklines
\put(5412,1668){\ellipse{2400}{450}}
\path(387,3918)(12,3843)(387,3693)
\path(11037,3993)(11412,3918)(11037,3768)
\path(5187,3468)(4962,3768)(4962,3618)
\path(4962,3768)(5112,3693)
\path(3087,3243)	(3098.816,3278.782)
	(3117.754,3321.233)
	(3162.000,3393.000)

\path(3162,3393)	(3206.598,3434.727)
	(3266.625,3479.332)
	(3330.590,3518.272)
	(3387.000,3543.000)

\path(3387,3543)	(3454.619,3553.418)
	(3495.437,3554.848)
	(3538.174,3554.385)
	(3580.764,3552.499)
	(3621.142,3549.660)
	(3687.000,3543.000)

\path(3687,3543)	(3754.152,3529.495)
	(3794.917,3519.001)
	(3837.686,3507.307)
	(3880.370,3495.388)
	(3920.878,3484.216)
	(3987.000,3468.000)

\path(3987,3468)	(4038.735,3458.441)
	(4102.160,3447.664)
	(4173.414,3436.303)
	(4210.771,3430.603)
	(4248.637,3424.995)
	(4286.531,3419.559)
	(4323.969,3414.375)
	(4395.547,3405.079)
	(4459.511,3397.742)
	(4512.000,3393.000)

\path(4512,3393)	(4563.935,3389.351)
	(4627.238,3385.288)
	(4698.182,3381.539)
	(4773.038,3378.832)
	(4810.767,3378.098)
	(4848.076,3377.897)
	(4919.568,3379.461)
	(4983.786,3384.252)
	(5037.000,3393.000)

\path(5037,3393)	(5106.783,3420.911)
	(5147.041,3442.278)
	(5188.710,3465.911)
	(5230.165,3489.806)
	(5269.782,3511.956)
	(5337.000,3543.000)

\path(5337,3543)	(5381.191,3555.343)
	(5435.443,3567.426)
	(5496.448,3578.954)
	(5560.897,3589.635)
	(5625.485,3599.174)
	(5686.903,3607.277)
	(5741.844,3613.650)
	(5787.000,3618.000)

\path(5787,3618)	(5853.188,3619.597)
	(5893.993,3618.998)
	(5937.000,3618.000)
	(5980.007,3617.002)
	(6020.812,3616.403)
	(6087.000,3618.000)

\path(6087,3618)	(6132.157,3622.538)
	(6187.098,3629.370)
	(6248.517,3638.039)
	(6313.106,3648.090)
	(6377.557,3659.067)
	(6438.561,3670.515)
	(6492.811,3681.978)
	(6537.000,3693.000)

\path(6537,3693)	(6574.744,3706.325)
	(6619.216,3725.885)
	(6668.356,3749.336)
	(6720.105,3774.338)
	(6772.403,3798.547)
	(6823.192,3819.622)
	(6870.410,3835.220)
	(6912.000,3843.000)

\path(6912,3843)	(6984.731,3842.803)
	(7025.914,3838.991)
	(7069.576,3833.055)
	(7115.168,3825.264)
	(7162.142,3815.881)
	(7209.947,3805.174)
	(7258.035,3793.406)
	(7305.856,3780.845)
	(7352.861,3767.755)
	(7398.501,3754.403)
	(7442.227,3741.054)
	(7483.488,3727.973)
	(7521.737,3715.427)
	(7587.000,3693.000)

\path(7587,3693)	(7655.814,3665.060)
	(7697.284,3646.247)
	(7740.521,3625.515)
	(7783.318,3603.843)
	(7823.468,3582.213)
	(7887.000,3543.000)

\path(7887,3543)	(7926.615,3507.615)
	(7962.000,3468.000)

\path(7962,3468)	(8011.481,3401.602)
	(8032.280,3360.104)
	(8037.000,3318.000)

\path(8037,3318)	(8014.493,3269.365)
	(7972.935,3227.887)
	(7925.909,3193.966)
	(7887.000,3168.000)

\path(7887,3168)	(7814.858,3125.363)
	(7772.641,3105.419)
	(7737.000,3093.000)

\path(7737,3093)	(7686.696,3086.991)
	(7623.938,3086.381)
	(7561.461,3089.081)
	(7512.000,3093.000)

\path(7512,3093)	(7474.161,3097.689)
	(7428.261,3104.858)
	(7377.040,3113.936)
	(7323.236,3124.354)
	(7269.591,3135.539)
	(7218.843,3146.923)
	(7173.733,3157.933)
	(7137.000,3168.000)

\path(7137,3168)	(7099.176,3181.303)
	(7054.081,3200.055)
	(7004.058,3222.309)
	(6951.450,3246.120)
	(6898.598,3269.541)
	(6847.844,3290.625)
	(6801.531,3307.427)
	(6762.000,3318.000)

\path(6762,3318)	(6695.539,3326.042)
	(6654.974,3328.337)
	(6612.293,3329.329)
	(6569.575,3328.905)
	(6528.900,3326.951)
	(6462.000,3318.000)

\path(6462,3318)	(6426.976,3304.232)
	(6385.789,3282.041)
	(6312.000,3243.000)

\path(6312,3243)	(6246.232,3224.973)
	(6205.216,3216.156)
	(6161.884,3207.270)
	(6118.567,3198.162)
	(6077.594,3188.682)
	(6012.000,3168.000)

\path(6012,3168)	(5960.438,3137.303)
	(5899.500,3093.000)
	(5838.562,3048.697)
	(5787.000,3018.000)

\path(5787,3018)	(5721.654,2996.007)
	(5680.898,2985.380)
	(5637.776,2975.235)
	(5594.558,2965.748)
	(5553.510,2957.095)
	(5487.000,2943.000)

\path(5487,2943)	(5450.063,2934.641)
	(5404.857,2924.433)
	(5354.097,2913.154)
	(5300.497,2901.585)
	(5246.773,2890.505)
	(5195.639,2880.695)
	(5149.810,2872.933)
	(5112.000,2868.000)

\path(5112,2868)	(5062.456,2864.421)
	(4999.841,2862.289)
	(4937.055,2863.012)
	(4887.000,2868.000)

\path(4887,2868)	(4847.426,2878.492)
	(4801.007,2895.095)
	(4750.123,2915.935)
	(4697.152,2939.134)
	(4644.476,2962.816)
	(4594.472,2985.105)
	(4549.520,3004.126)
	(4512.000,3018.000)

\path(4512,3018)	(4446.700,3037.680)
	(4406.227,3048.996)
	(4363.417,3060.386)
	(4320.431,3071.165)
	(4279.426,3080.650)
	(4212.000,3093.000)

\path(4212,3093)	(4145.634,3097.271)
	(4104.812,3098.019)
	(4061.820,3098.014)
	(4018.850,3097.381)
	(3978.096,3096.250)
	(3912.000,3093.000)

\path(3912,3093)	(3859.502,3088.452)
	(3795.517,3081.586)
	(3723.913,3072.872)
	(3686.463,3067.969)
	(3648.559,3062.779)
	(3610.684,3057.362)
	(3573.322,3051.775)
	(3502.071,3040.332)
	(3438.674,3028.917)
	(3387.000,3018.000)

\path(3387,3018)	(3341.044,2995.190)
	(3285.525,2960.329)
	(3224.494,2935.553)
	(3162.000,2943.000)

\path(3162,2943)	(3116.731,2987.826)
	(3096.908,3050.347)
	(3090.880,3115.446)
	(3087.000,3168.000)

\path(3087,3168)	(3084.375,3205.384)
	(3087.000,3243.000)

\path(4212,1668)	(4215.527,1723.964)
	(4216.703,1765.275)
	(4212.000,1818.000)

\path(4212,1818)	(4186.260,1889.002)
	(4166.793,1930.327)
	(4144.879,1972.961)
	(4121.979,2014.975)
	(4099.558,2054.440)
	(4062.000,2118.000)

\path(4062,2118)	(4028.683,2172.710)
	(3987.381,2238.913)
	(3940.148,2312.420)
	(3914.951,2350.605)
	(3889.043,2389.046)
	(3862.679,2427.220)
	(3836.119,2464.604)
	(3783.433,2534.907)
	(3733.042,2595.768)
	(3687.000,2643.000)

\path(3687,2643)	(3638.550,2681.385)
	(3576.904,2722.430)
	(3506.206,2764.678)
	(3468.758,2785.798)
	(3430.601,2806.673)
	(3392.254,2827.120)
	(3354.233,2846.958)
	(3281.248,2884.079)
	(3215.788,2916.578)
	(3162.000,2943.000)

\path(3162,2943)	(3121.397,2962.936)
	(3075.186,2985.512)
	(3024.122,3010.261)
	(2968.958,3036.719)
	(2910.451,3064.420)
	(2849.354,3092.897)
	(2786.422,3121.687)
	(2722.410,3150.323)
	(2658.073,3178.339)
	(2594.165,3205.270)
	(2531.441,3230.650)
	(2470.657,3254.015)
	(2412.565,3274.897)
	(2357.922,3292.833)
	(2307.482,3307.356)
	(2262.000,3318.000)

\path(2262,3318)	(2194.462,3324.607)
	(2155.269,3324.590)
	(2113.203,3322.895)
	(2068.821,3319.924)
	(2022.680,3316.080)
	(1975.336,3311.763)
	(1927.346,3307.376)
	(1879.268,3303.322)
	(1831.657,3300.001)
	(1785.071,3297.816)
	(1740.066,3297.169)
	(1697.200,3298.462)
	(1657.029,3302.097)
	(1587.000,3318.000)

\path(1587,3318)	(1531.930,3346.629)
	(1474.930,3390.265)
	(1416.796,3444.012)
	(1358.325,3502.972)
	(1300.313,3562.249)
	(1243.557,3616.943)
	(1188.854,3662.160)
	(1137.000,3693.000)

\path(1137,3693)	(1099.465,3708.448)
	(1059.483,3722.830)
	(1016.641,3736.201)
	(970.527,3748.615)
	(920.729,3760.129)
	(866.836,3770.796)
	(808.435,3780.671)
	(745.114,3789.810)
	(676.461,3798.267)
	(602.065,3806.098)
	(562.584,3809.796)
	(521.513,3813.357)
	(478.800,3816.790)
	(434.394,3820.100)
	(388.243,3823.294)
	(340.295,3826.380)
	(290.499,3829.364)
	(238.804,3832.254)
	(185.158,3835.055)
	(129.510,3837.775)
	(71.808,3840.421)
	(12.000,3843.000)

\path(6612,1668)	(6609.621,1724.223)
	(6608.828,1765.620)
	(6612.000,1818.000)

\path(6612,1818)	(6623.229,1869.479)
	(6641.505,1932.776)
	(6663.778,1994.936)
	(6687.000,2043.000)

\path(6687,2043)	(6730.342,2099.964)
	(6790.901,2164.418)
	(6855.760,2224.412)
	(6912.000,2268.000)

\path(6912,2268)	(6948.018,2286.949)
	(6993.502,2304.412)
	(7045.269,2321.102)
	(7100.138,2337.735)
	(7154.927,2355.026)
	(7206.455,2373.691)
	(7251.540,2394.444)
	(7287.000,2418.000)

\path(7287,2418)	(7318.636,2451.724)
	(7349.419,2496.046)
	(7379.275,2547.546)
	(7408.129,2602.804)
	(7435.903,2658.398)
	(7462.524,2710.909)
	(7487.915,2756.917)
	(7512.000,2793.000)

\path(7512,2793)	(7552.556,2842.074)
	(7603.976,2899.830)
	(7663.069,2962.952)
	(7726.646,3028.125)
	(7791.518,3092.032)
	(7854.494,3151.358)
	(7912.384,3202.785)
	(7962.000,3243.000)

\path(7962,3243)	(8025.134,3287.863)
	(8062.761,3312.902)
	(8103.650,3339.188)
	(8147.221,3366.355)
	(8192.891,3394.033)
	(8240.079,3421.855)
	(8288.201,3449.453)
	(8336.677,3476.458)
	(8384.924,3502.503)
	(8432.359,3527.220)
	(8478.402,3550.240)
	(8522.469,3571.197)
	(8563.979,3589.721)
	(8637.000,3618.000)

\path(8637,3618)	(8673.645,3629.516)
	(8713.159,3640.671)
	(8755.352,3651.471)
	(8800.033,3661.925)
	(8847.011,3672.038)
	(8896.096,3681.820)
	(8947.096,3691.276)
	(8999.821,3700.413)
	(9054.081,3709.240)
	(9109.685,3717.763)
	(9166.442,3725.990)
	(9224.161,3733.927)
	(9282.652,3741.582)
	(9341.724,3748.962)
	(9401.187,3756.074)
	(9460.849,3762.926)
	(9520.520,3769.525)
	(9580.009,3775.877)
	(9639.126,3781.991)
	(9697.680,3787.872)
	(9755.480,3793.530)
	(9812.336,3798.969)
	(9868.057,3804.199)
	(9922.452,3809.226)
	(9975.330,3814.057)
	(10026.501,3818.700)
	(10075.775,3823.161)
	(10122.960,3827.448)
	(10167.865,3831.569)
	(10210.301,3835.529)
	(10250.076,3839.337)
	(10287.000,3843.000)

\path(10287,3843)	(10324.102,3846.645)
	(10363.710,3850.301)
	(10406.235,3853.995)
	(10452.090,3857.755)
	(10501.685,3861.608)
	(10555.435,3865.582)
	(10613.749,3869.704)
	(10677.041,3874.001)
	(10745.723,3878.502)
	(10820.205,3883.233)
	(10859.751,3885.694)
	(10900.901,3888.222)
	(10943.708,3890.822)
	(10988.222,3893.497)
	(11034.496,3896.250)
	(11082.581,3899.085)
	(11132.528,3902.004)
	(11184.389,3905.013)
	(11238.215,3908.113)
	(11294.058,3911.309)
	(11351.969,3914.603)
	(11412.000,3918.000)

\path(4212,1668)	(4215.527,1612.036)
	(4216.703,1570.725)
	(4212.000,1518.000)

\path(4212,1518)	(4199.051,1482.533)
	(4178.224,1440.608)
	(4137.000,1368.000)

\path(4137,1368)	(4106.015,1317.039)
	(4065.169,1253.501)
	(4022.738,1190.963)
	(3987.000,1143.000)

\path(3987,1143)	(3955.290,1108.417)
	(3913.189,1066.635)
	(3870.494,1025.536)
	(3837.000,993.000)

\path(3837,993)	(3802.951,961.301)
	(3759.330,921.596)
	(3717.044,880.094)
	(3687.000,843.000)

\path(3687,843)	(3662.524,795.740)
	(3638.160,734.921)
	(3619.466,671.892)
	(3612.000,618.000)

\path(3612,618)	(3617.825,564.022)
	(3633.784,500.850)
	(3657.600,440.002)
	(3687.000,393.000)

\path(3687,393)	(3748.892,344.873)
	(3789.242,323.268)
	(3832.594,303.409)
	(3876.496,285.373)
	(3918.498,269.240)
	(3987.000,243.000)

\path(3987,243)	(4030.786,225.574)
	(4084.644,205.316)
	(4145.319,183.451)
	(4209.555,161.209)
	(4274.096,139.815)
	(4335.688,120.497)
	(4391.074,104.483)
	(4437.000,93.000)

\path(4437,93)	(4503.615,80.474)
	(4542.876,74.327)
	(4585.275,68.302)
	(4630.189,62.433)
	(4676.994,56.756)
	(4725.068,51.303)
	(4773.788,46.110)
	(4822.529,41.210)
	(4870.670,36.638)
	(4917.587,32.428)
	(4962.656,28.613)
	(5005.256,25.229)
	(5044.761,22.310)
	(5112.000,18.000)

\path(5112,18)	(5160.663,15.667)
	(5216.111,13.645)
	(5277.374,11.934)
	(5343.484,10.534)
	(5413.470,9.445)
	(5486.365,8.667)
	(5561.198,8.200)
	(5599.038,8.083)
	(5637.000,8.044)
	(5674.962,8.082)
	(5712.802,8.199)
	(5787.635,8.665)
	(5860.530,9.443)
	(5930.516,10.531)
	(5996.626,11.931)
	(6057.889,13.643)
	(6113.337,15.666)
	(6162.000,18.000)

\path(6162,18)	(6229.022,20.566)
	(6268.292,21.586)
	(6310.601,22.635)
	(6355.350,23.866)
	(6401.938,25.429)
	(6449.765,27.480)
	(6498.233,30.169)
	(6546.739,33.650)
	(6594.686,38.075)
	(6641.472,43.596)
	(6686.497,50.368)
	(6729.163,58.541)
	(6768.869,68.269)
	(6837.000,93.000)

\path(6837,93)	(6897.325,134.677)
	(6960.945,196.136)
	(7018.842,262.277)
	(7062.000,318.000)

\path(7062,318)	(7104.873,378.357)
	(7129.707,415.728)
	(7154.276,456.011)
	(7176.686,497.793)
	(7195.041,539.660)
	(7207.444,580.200)
	(7212.000,618.000)

\path(7212,618)	(7203.521,657.168)
	(7183.140,698.992)
	(7137.000,768.000)

\path(7137,768)	(7106.434,804.014)
	(7064.846,846.214)
	(7021.835,886.807)
	(6987.000,918.000)

\path(6987,918)	(6920.529,966.097)
	(6878.162,993.658)
	(6833.441,1022.775)
	(6789.165,1052.858)
	(6748.133,1083.316)
	(6687.000,1143.000)

\path(6687,1143)	(6661.745,1190.943)
	(6639.698,1253.119)
	(6622.552,1316.485)
	(6612.000,1368.000)

\path(6612,1368)	(6607.236,1412.580)
	(6605.647,1472.753)
	(6606.045,1510.883)
	(6607.236,1555.549)
	(6609.221,1607.628)
	(6612.000,1668.000)

\put(4812,3768){\makebox(0,0)[lb]{\smash{{{\SetFigFont{12}{14.4}{rm}$\vec{\nu}$}}}}}
\put(10512,3293){\makebox(0,0)[lb]{\smash{{{\SetFigFont{12}{14.4}{rm}$(M^3,g_0)$}}}}}
\put(4212,1043){\makebox(0,0)[lb]{\smash{{{\SetFigFont{12}{14.4}{rm}$\Sigma(0) = \Sigma_0$}}}}}
\put(2337,2568){\makebox(0,0)[lb]{\smash{{{\SetFigFont{12}{14.4}{rm}$\Sigma(t)$}}}}}
\end{picture}